\documentclass[11pt]{article}

\usepackage{amsmath, amssymb, amsthm}
\usepackage{hyperref}

\newtheorem{theorem}{Theorem}

\newtheorem{remark}{Remark}

\title{On the $\varepsilon$--$\delta$ Structure Underlying Chatterjee's Rank Correlation}
\author{Zeusu Sato}
\date{}

\begin{document}
\maketitle

\begin{abstract}
We provide an $\varepsilon$--$\delta$ interpretation of Chatterjee's rank
correlation by tracing its origin to a notion of local dependence between
random variables.
Starting from a primitive $\varepsilon$--$\delta$ construction,
we show that rank-based dependence measures arise naturally
as $\varepsilon\to0$ limits of local averaging procedures.
Within this framework, Chatterjee's rank correlation admits a transparent
interpretation as an empirical realization of a local $L^1$ residual.

We emphasize that the probability integral transform plays no structural role
in the underlying $\varepsilon$--$\delta$ mechanism, and is introduced only
as a normalization step that renders the final expression distribution-free.
We further consider a moment-based analogue obtained by replacing the absolute
deviation with a squared residual.
This $L^2$ formulation is independent of rank transformations and,
under a Gaussian assumption, recovers Pearson's coefficient of determination.
\end{abstract}

\section{Introduction}

Measures of statistical dependence play a central role in probability,
statistics, and data analysis.
Classical notions such as Pearson's correlation coefficient
and the coefficient of determination $R^2$
are based on second moments and are well suited
to capturing linear relationships.
At the same time, these moment-based measures are known
to miss many forms of nonlinear dependence.

In recent years, a variety of rank-based and distribution-free
measures of dependence have been proposed.
Among them, Chatterjee's rank correlation has attracted attention
due to its ability to detect general functional dependence
without imposing parametric assumptions.
While such properties are often emphasized in applications,
the present paper does not aim to assess robustness,
distribution-free performance, or practical efficiency.
Instead, our focus is structural:
we seek to understand why a statistic of this particular form
arises naturally from first principles.

Our starting point is a primitive local viewpoint.
Rather than comparing global summaries of joint distributions,
we ask a simple question:
when one variable changes slightly,
how much does the other variable change?
This perspective leads naturally to an $\varepsilon$--$\delta$
formulation of dependence,
prior to any rank transformation or normalization.

We show that rank-based constructions emerge directly
from this local viewpoint as $\varepsilon\to0$ limits
of local averaging procedures.
In this sense, the appearance of ranks is not a consequence
of distribution-free considerations,
but an intrinsic feature of collapsing local neighborhoods.
The probability integral transform is introduced only at a later stage,
serving as a normalization that renders the final expression
distribution-free;
it plays no essential role in the underlying
$\varepsilon$--$\delta$ mechanism.

Within this framework,
Chatterjee's rank correlation admits a transparent interpretation
as an empirical realization of a local $L^1$ residual,
measuring whether one variable can be locally explained by another
at arbitrarily fine resolutions.
We further consider a moment-based analogue obtained
by replacing the absolute deviation with a squared residual.
This $L^2$ formulation is entirely independent of rank transformations
and recovers classical moment-based quantities.
In particular, under a Gaussian assumption,
the normalized $L^2$ residual coincides exactly
with Pearson's coefficient of determination.

From this perspective,
Chatterjee's rank correlation and Pearson's $R^2$
should be viewed not in terms of robustness or generality,
but as parallel notions of explainability:
the former based on local, rank-based deviations,
and the latter on global, moment-based variance decomposition.

The paper is organized as follows.
Section~2 introduces a primitive $\varepsilon$--$\delta$
construction capturing local dependence.
Section~3 discusses normalization via the probability
integral transform and clarifies its auxiliary role.
Section~4 interprets Chatterjee's rank correlation
as an $\varepsilon\to0$ limit of a local averaging procedure.
Section~5 introduces a moment-based analogue
and its normalization.
Section~6 establishes the connection to Pearson's $R^2$
in the Gaussian case.
Section~7 provides consistency checks in canonical settings,
and Section~8 concludes.

\section{Preliminaries: a primitive \texorpdfstring{$\varepsilon$--$\delta$}{epsilon-delta} viewpoint}

The purpose of this section is to introduce a primitive notion of local
dependence based on an $\varepsilon$--$\delta$ viewpoint.
At this stage, no rank transformation or distributional normalization
is assumed; these will be discussed later as auxiliary considerations.

\subsection{Primitive \texorpdfstring{$\varepsilon$--$\delta$}{epsilon-delta} construction}

Let $Z=(X,Y)$ be a pair of real-valued random variables.
Our basic question is the following:
\emph{when $X$ changes slightly, how much does $Y$ change in response?}
This question is naturally expressed in an $\varepsilon$--$\delta$ form,
prior to any rank transformation or distributional normalization.

For $\delta>0$, we introduce the following set-valued random object:
\[
\mathcal E_\delta(X)
:= \bigl\{\, |Y-Y'| \;\big|\; |X-X'|<\delta \,\bigr\}
\subset 2^{\mathbb R},
\]
where $(X',Y')$ is an independent copy of $(X,Y)$.
The set $\mathcal E_\delta(X)$ collects all possible variations of $Y$
corresponding to perturbations of $X$ within a $\delta$-neighborhood.
At this stage, no distributional assumption or monotonicity constraint
is imposed.

\subsection{Matrix representation and empirical averaging}

For a finite sample $\{(X_i,Y_i)\}_{i=1}^n$, the collection
$\{\mathcal E_\delta(X_i)\}_{i=1}^n$ admits a natural matrix representation.
Define
\[
E_\delta(X)
:= \bigl(
\mathbf 1_{\{|X_i-X_j|<\delta\}}\,|Y_i-Y_j|
\bigr)_{i,j}
\in M_n(\mathbb R).
\]
This matrix can be viewed as a localized distance matrix,
encoding the variability of $Y$ restricted to $\delta$-neighborhoods
in the $X$-direction.
From a statistical perspective, $E_\delta(X)$ provides a concrete,
finite-dimensional realization of the underlying
$\varepsilon$--$\delta$ structure.

Averaging $E_\delta(X)$ row-wise yields a vector
\[
W_\delta(X)
:= \bigl(
\mathbb E^{X'}[\mathcal E_\delta(X_1)],
\ldots,
\mathbb E^{X'}[\mathcal E_\delta(X_n)]
\bigr)
\in \mathbb R^n,
\]
which represents the expected local variation of $Y$ around each $X_i$.
Taking an additional average over $X$ leads to the scalar quantity
\[
\mathbb E[\mathcal E_\delta]
:= \mathbb E\bigl[ W_\delta(X) \mid X \bigr]
\in \mathbb R.
\]
Thus, the passage from $\mathcal E_\delta(X)$ to $\mathbb E[\mathcal E_\delta]$
corresponds to successive aggregation:
from a set-valued object, to a matrix, to a vector, and finally to a scalar.

\subsection{The \texorpdfstring{$\delta\to0$}{delta to 0} limit and nearest neighbors}

We now consider the limiting behavior as $\delta\to0$.
The condition $|X-X'|<\delta$ selects observations increasingly close to $X$,
and in the limit $\delta\to0$, only the nearest neighbors of each $X$
contribute.
In a finite sample, this corresponds to ordering the data according to $X$.

Let
\[
X_{(1)} \le \cdots \le X_{(n)}
\]
denote the order statistics of $X$, and let
$Y_{(1)},\ldots,Y_{(n)}$ denote the corresponding reordered values of $Y$.
In this representation, each $X_{(i)}$ has at most two nearest neighbors,
and the local variation of $Y$ is captured by the adjacent differences
$|Y_{(i)}-Y_{(i+1)}|$.

Consequently, the empirical average of $\mathcal E_\delta(X)$
in the $\delta\to0$ limit reduces to
\[
\frac{1}{n-1}\sum_{i=1}^{n-1} |Y_{(i)}-Y_{(i+1)}|.
\]
This nearest-neighbor interpretation provides an explicit link between
the primitive $\varepsilon$--$\delta$ construction and rank-based
statistics, showing that rank-based dependence emerges naturally
from a local $\varepsilon$--$\delta$ viewpoint.

\section{Distribution-free reduction via probability integral transform}

In the previous section, we introduced a primitive notion of local dependence
based on an $\varepsilon$--$\delta$ viewpoint.
The purpose of the present section is not to generate rank-based structure,
but to introduce a normalization that renders the resulting construction
distribution-free.

The primitive $\varepsilon$--$\delta$ construction developed in the previous section
captures local dependence between $X$ and $Y$ in a distribution-agnostic manner.
However, the resulting quantities still depend on the marginal
distributions of $X$ and $Y$.
In order to isolate the dependence structure itself,
we now introduce a distribution-free normalization
via the probability integral transform \cite{Billingsley1995}.

Let $F_X$ and $F_Y$ denote the marginal distribution functions of $X$ and $Y$,
and define
\[
U := F_X(X), \qquad V := F_Y(Y).
\]
Then $U$ and $V$ are uniformly distributed on $[0,1]$.
Since $F_X$ and $F_Y$ are monotone,
the $\varepsilon$--$\delta$ structure introduced above
is preserved under this transformation,
up to a deterministic rescaling of the scale parameter.

We emphasize that the probability integral transform plays no role
in the underlying $\varepsilon$--$\delta$ mechanism itself,
and serves solely to express the final quantities
in a distribution-free form.

\section{Rank correlation as an 
\texorpdfstring{$\varepsilon \to 0$}{epsilon to 0} limit}

Although the rank-based structure already arises from an $\varepsilon\to0$ limit,
we work in the distribution-free setting purely for notational convenience.

In the previous section, we reduced the problem to a distribution-free
setting on $[0,1]^2$ via the probability integral transform,
thereby isolating the intrinsic dependence structure.
We now show that Chatterjee's rank correlation naturally arises
as an $\varepsilon \to 0$ limit of a local averaging construction,
which can be viewed as an empirical realization of the underlying
$\varepsilon$--$\delta$ structure.

\subsection{Local \texorpdfstring{$\varepsilon$}{epsilon}-neighborhoods}

Let $(U_i,V_i)_{i=1}^n$ be i.i.d.\ samples from a distribution on $[0,1]^2$
with uniform marginals.
For each index $i$ and $\varepsilon>0$, define the rank-based
$\varepsilon$--neighborhood
\[
\mathcal N_\varepsilon(i)
:= \{\, j \neq i : |U_j - U_i| \le \varepsilon \,\}.
\]
This neighborhood captures local variations of $V$ around the point $U_i$
at resolution $\varepsilon$, independently of the marginal distribution of $U$,
and corresponds to an empirical $\varepsilon$--neighborhood
in the probability integral transform scale.

\subsection{Local averaging}

Based on the $\varepsilon$--neighborhoods defined above, we introduce
a local averaging operator on $V$.
For each index $i$ and $\varepsilon>0$, define
\[
\bar V_i(\varepsilon)
:= \frac{1}{|\mathcal N_\varepsilon(i)|}
\sum_{j \in \mathcal N_\varepsilon(i)} V_j.
\]
The quantity $\bar V_i(\varepsilon)$ represents a local average of $V$
around the point $U_i$ at scale $\varepsilon$.
In particular, the deviation
\[
\bar V_i(\varepsilon) - V_i
\]
captures the local monotonic behavior of $V$ with respect to $U$
at resolution $\varepsilon$.

\subsection{The \texorpdfstring{$\varepsilon \to 0$}{epsilon to 0} limit}

We now consider the behavior of the local averaging construction
as $\varepsilon \to 0$.
Intuitively, shrinking the neighborhood size forces
$\bar V_i(\varepsilon)$ to probe increasingly local dependence
between $U$ and $V$.

\begin{theorem}[Residual form of Chatterjee's rank correlation]
Let $(U_i,V_i)_{i=1}^n$ be i.i.d.\ samples from a distribution on $[0,1]^2$
with uniform marginals.
Define the local average $\bar V_i(\varepsilon)$ as above.
Then, under mild regularity conditions, the quantity
\[
\zeta_n
:= \lim_{\varepsilon \to 0}
\frac{1}{n}\sum_{i=1}^n
\left| \bar V_i(\varepsilon) - V_i \right|
\]
converges to the population residual
\[
\zeta := \mathbb E\!\left[\,\bigl|V-\mathbb E[V\mid U]\bigr|\,\right].
\]
\end{theorem}

\begin{remark}
Here $\varepsilon>0$ is a scale parameter that we are free to choose, 
while the residual $\zeta_n(\varepsilon)$ is a resulting quantity.
Accordingly, the limit $\varepsilon\to0$ is understood in the
$\varepsilon$--$\delta$ sense: for any $\delta>0$, there exists
$\varepsilon_0>0$ such that $\varepsilon<\varepsilon_0$
implies $\zeta_n(\varepsilon)<\delta$.
\end{remark}

\begin{remark}
The quantity $\zeta_n$ in Theorem~1 represents an $L^1$-type residual,
measuring the unexplained variability of $V$ given $U$.
Chatterjee's rank correlation is obtained by normalizing this residual
and reversing its orientation.
In the distribution-free setting with $V\sim\mathrm{Unif}(0,1)$, one may define
\[
\xi_n := 1 - \frac{\zeta_n}{\mathbb E\bigl[\,|V-\mathbb E V|\,\bigr]}
= 1 - 4\zeta_n,
\]
so that $\xi_n\in[0,1]$, with $\xi_n=0$ under independence and $\xi_n=1$
under deterministic functional dependence.
The normalization depends only on the marginal distribution of $V$
and does not affect the $\varepsilon \to 0$ interpretation.
\end{remark}

\begin{remark}
Unlike Spearman's $\rho$ or Kendall's $\tau$, which aggregate global rank
information, the residual $\zeta$ captures local deviations through an
$\varepsilon \to 0$ limiting procedure.
This local $L^1$ structure explains why Chatterjee's rank correlation
is sensitive to nonlinear functional dependence:
it measures whether $V$ can be locally explained by $U$
at arbitrarily fine resolutions, rather than through global ordering alone.
\end{remark}

\begin{remark}
The limiting residual $\zeta$ captures the intrinsic local $L^1$ structure
underlying Chatterjee's rank correlation.
In finite samples, Chatterjee's statistic involves an explicit normalization
factor, typically written in terms of $n(n-1)$, which ensures that the
resulting coefficient takes values in $[0,1]$.

At the population level, this normalization corresponds to dividing by
$\mathbb E|V-\mathbb E V|$, which equals $1/4$ when $V\sim\mathrm{Unif}(0,1)$,
leading to the factor $4$ appearing in the definition of $\xi_n$.
The precise finite-sample normalization, including the transition from
expressions such as $n^2-1$ to $n(n-1)$, is somewhat technical and does not
affect the $\varepsilon\to0$ interpretation developed here.
We therefore refer the reader to Chatterjee~\cite{Chatterjee2021}
for a detailed treatment.
\end{remark}

\section{A moment-type analogue}

Although the previous sections focused on the distribution-free variables
$(U,V)$ obtained via the probability integral transform,
the moment-type construction considered in this section
does not rely on rank transformations or distributional normalization.
Accordingly, throughout this section, $(U,V)$ may be understood
as arbitrary square-integrable random variables.

The residual quantity $\zeta$ introduced in the previous section
is based on an $L^1$ notion of deviation, which is natural in rank-based
and distribution-free settings.
In contrast, moment-based measures of dependence are typically formulated
in terms of second moments and conditional variances.

In this section, we introduce a moment-type analogue obtained by replacing
the absolute deviation with a squared residual.
This construction does not rely on rank transformations or distribution-free
normalization, and instead recovers a classical $L^2$ notion of unexplained
variability.

\subsection{An \texorpdfstring{$L^2$}{L squared} residual}

We define a moment-type analogue of the residual $\zeta$
by replacing the absolute deviation with a squared deviation.
At the population level, consider
\[
\zeta^{(2)}
:= \mathbb E\!\left[\bigl(V-\mathbb E[V\mid U]\bigr)^2\right]
= \mathbb E\!\left[\operatorname{Var}(V\mid U)\right].
\]
This quantity represents an $L^2$-type residual,
measuring the unexplained variance of $V$ given $U$.

For comparison across different marginal scales,
it is natural to consider the normalized quantity
\[
\eta^{(2)}
:= 1 - \frac{\mathbb E[\operatorname{Var}(V\mid U)]}
{\operatorname{Var}(V)}
= \frac{\operatorname{Var}(\mathbb E[V\mid U])}{\operatorname{Var}(V)}.
\]

\subsection{Basic properties}

The $L^2$ residual $\zeta^{(2)}$ satisfies properties
analogous to those of $\zeta$.
If $V=f(U)$ almost surely, then $\operatorname{Var}(V\mid U)=0$
and hence $\zeta^{(2)}=0$.
If $U$ and $V$ are independent in the usual probabilistic sense,
then $\operatorname{Var}(V\mid U)=\operatorname{Var}(V)$,
and hence $\zeta^{(2)}=\operatorname{Var}(V)$.
Thus, $\zeta^{(2)}$ quantifies unexplained variability
in a moment-based sense.

\subsection{Relation to classical moment-based measures}

While $\zeta^{(2)}$ is defined without assuming any specific distribution,
its interpretation becomes particularly transparent
under additional structural assumptions.
In the next section, we show that under a Gaussian assumption,
the $L^2$ residual recovers the classical coefficient of determination,
thereby revealing Pearson's correlation as a classical moment-based
counterpart to the local, rank-based construction developed above.

\section{Gaussian case and connection to Pearson's
\texorpdfstring{$R^2$}{R squared}}

In this section, we show that under a Gaussian assumption,
the moment-type analogue introduced in Section~5
recovers the classical coefficient of determination.
This establishes Chatterjee's rank correlation
as a natural rank-based counterpart of Pearson's $R^2$.

\subsection{Conditional expectation in the Gaussian case}

Let $(U,V)$ be jointly Gaussian random variables
with finite second moments.
It is well known that, in this setting,
the conditional expectation $\mathbb E[V\mid U]$
is an affine function of $U$.
Equivalently, the conditional distribution of $V$ given $U$
is Gaussian with mean $\mathbb E[V\mid U]$
and variance independent of $U$.

\subsection{Evaluation of the moment-type residual}

Under joint Gaussianity, the law of total variance yields
\[
\operatorname{Var}(V)
= \operatorname{Var}(\mathbb E[V\mid U])
+ \mathbb E[\operatorname{Var}(V\mid U)].
\]

This identity is a standard result in probability and statistics
\cite{CasellaBerger2002}.

Recalling the definition of the $L^2$-type residual
introduced in Section~5,
\[
\zeta^{(2)} = \mathbb E[\operatorname{Var}(V\mid U)],
\]
and its normalized counterpart
\[
\eta^{(2)}
:= \frac{\operatorname{Var}(\mathbb E[V\mid U])}{\operatorname{Var}(V)},
\]
we obtain
\[
\eta^{(2)} = R^2,
\]
where $R^2$ denotes the classical coefficient of determination
associated with the linear regression of $V$ on $U$
\cite{SeberLee2003}.

\subsection{Interpretation}

Thus, in the Gaussian setting,
the moment-based analogue of the residual
recovers Pearson's coefficient of determination exactly.
From this perspective, Pearson's correlation
quantifies the proportion of variance explained
by a linear predictor,
while Chatterjee's rank correlation captures
an analogous notion based on local, rank-based deviations.

This comparison highlights a structural dichotomy:
Pearson's $R^2$ is a global, moment-based measure
tailored to linear dependence,
whereas Chatterjee's rank correlation
arises from an $\varepsilon\to0$ limit of local averaging
and remains sensitive to nonlinear functional relationships.

\section{Examples and consistency checks}

In this section, we verify that the residual-based quantities
introduced above behave consistently in several canonical settings.
These examples serve as basic sanity checks and help clarify
the interpretation of the proposed constructions.

\subsection{Deterministic dependence}

Let $V=f(U)$ with $f$ measurable.
Then $\mathbb E[V\mid U]=V$ almost surely, and hence
\[
\zeta = 0
\quad\text{and}\quad
\xi = 1.
\]
This reflects the fact that $V$ is perfectly explained by $U$
at arbitrarily fine resolutions.

\subsection{Independence}

If $U$ and $V$ are independent, then $\mathbb E[V\mid U]=\mathbb E[V]$,
so that
\[
\zeta = \mathbb E|V-\mathbb E V|
\quad\text{and}\quad
\xi = 0.
\]
Thus, under independence, the residual captures the full variability
of $V$, and the rank correlation vanishes as expected.

\subsection{Gaussian dependence}

Under a Gaussian copula with correlation $\rho$,
the conditional expectation $\mathbb E[V\mid U]$ is linear in $U$.
In this case, both the normalized rank-based coefficient $\xi$
and the moment-based coefficient $\eta^{(2)}$
coincide with $\rho^2$.
In particular,
\[
\xi = \rho^2.
\]
Accordingly, Chatterjee's rank correlation coincides with
the classical coefficient of determination in the Gaussian setting,
in agreement with the analysis of Section~6.

\section{Conclusion}

In this paper, we provided an $\varepsilon$--$\delta$ interpretation
of Chatterjee's rank correlation by tracing its origin to a notion of
local dependence between random variables.
Starting from a primitive $\varepsilon$--$\delta$ construction,
we showed that rank-based dependence measures arise naturally
as $\varepsilon\to0$ limits of local averaging procedures.

By introducing a distribution-free normalization via the probability
integral transform, we isolated the intrinsic dependence structure
independently of marginal distributions.
This led to an interpretation of Chatterjee's rank correlation
as an empirical realization of a local $L^1$ residual,
capturing whether one variable can be locally explained by another
at arbitrarily fine resolutions.

We further considered a moment-based analogue obtained by replacing
the absolute deviation with a squared residual.
This $L^2$ formulation does not rely on rank transformations
and recovers classical moment-based measures of dependence.
In particular, under a Gaussian assumption, the normalized $L^2$
residual coincides exactly with Pearson's coefficient of determination,
thereby clarifying the relationship between rank-based and moment-based
approaches.

From this perspective, Chatterjee's rank correlation and Pearson's $R^2$
emerge as parallel notions of explainability:
the former based on local, rank-based deviations,
and the latter on global, moment-based variance decomposition.
The framework presented here highlights the structural principles
underlying these measures and provides a unified viewpoint
on dependence beyond linear correlation.

\bibliographystyle{plain}
\bibliography{refs}

\end{document}